\documentclass[12pt]{amsart}
\usepackage{amssymb}
\usepackage{mathrsfs}
\usepackage{cite}
\usepackage{algorithm}
\usepackage{algorithmic}
\usepackage{amscd}

{\theoremstyle{definition} 
}
\newtheorem{pro}{Proposition}
\newtheorem{lmm}{Lemma}
\newtheorem{thm}{Theorem}
\newtheorem{cor}{Corollary}

\begin{document}

\title{Finite subsets of projective space, and their ideals}
\author{Mathias Lederer}
\address{Fakult\"at f\"ur Mathematik, Universit\"at Bielefeld, P.O.Box 100131, D-33501 Bielefeld, Germany}
\email{mlederer@math.uni-bielefeld.de	}
\date{November 2007}

\begin{abstract}
  Let $\mathscr{A}$ be a finite set of closed rational points in projective space, 
  let $\mathscr{I}$ be the vanishing ideal of $\mathscr{A}$, 
  and let $\mathscr{D}(\mathscr{A})$
  be the set of exponents of those monomials which do not occur as leading monomials of elements of $\mathscr{I}$. 
  We show that the size of $\mathscr{A}$ equals the number of axes contained in $\mathscr{D}(\mathscr{A})$.
  Furthermore, we present an algorithm for the construction of the Gr\"obner basis of $\mathscr{I}(\mathscr{A})$, 
  hence also of $\mathscr{D}(\mathscr{A})$.
\end{abstract}

\maketitle

\section{Introduction}

The aim of this article is to carry over a well-known fact of affine geometry to projective geometry.
Let us first discuss the affine statement.

Consider the $n$-dimensional affine space $\mathbb{A}^n$ over a field $k$. 
Given an ideal $I\subset k[X]=k[X_{1},\ldots,X_{n}]$, and a term order $\leq$ on $k[X]$, 
we denote by $C(I)$ the set of all $\alpha\in\mathbb{N}^n$ which occur as exponents of leading terms of elements of $I$. 
That is, the monomial ideal $({\rm LT}(I))$ is given by $(X^\alpha;\alpha\in C(A))$. 
Furthermore, we consider the complement $D(I)=\mathbb{N}^n-C(I)$, 
which in the literature is called the set of {\it standard monomials} of $I$. 
The $k$-vector space $_{k}\langle X^\alpha;\alpha\in D(I)\rangle$ is isomorphic to the $k$-vector space $k[X]/I$, 
via the canonical map $X^\alpha\mapsto X^\alpha+I$. 
In this situation, the following statements holds: 
$D(I)$ is a finite set if, and only if, 
for all field extensions $L$ of $k$, the system of equations $f(a)=0$, for all $f\in I$, 
has only a finite number of solutions $a\in\mathbb{A}^n(L)$. 
If $D(I)$ is a finite set and $I$ is a radical ideal, 
then $\#D(I)$ equals the number of solutions of the system of equations $f(a)=0$, for all $f\in I$,
over $\overline{k}$, the algebraic closure of $k$. 

This follows from the same arguments as those provided in Chapter 5, \S3 in \cite{cox}. 
Note that the shape of $D(I)$ strongly depends on the term order $\leq$; 
however, the size of $D(I)$ depends only on $I$. 

Consider, conversely, a finite set $A$ of closed $k$-rational points of $\mathbb{A}^n$, 
and the ideal $I(A)$ consisting of those polynomials in $k[X_{1},\ldots,X_{n}]$ which vanish at all elements of $A$. 
We always use the shorthand notation $C(A)=C(I(A))$ and $D(A)=D(I(A))$. 
Then by the above, $D(A)$ is a finite set of size $\#A$. 
We are going to prove an analogue of this statement in a projective setting. 

We work in projective space $\mathbb{P}^n$ over the field $k$, 
in which we use the coordinates $X_{1},\ldots,X_{n+1}$. Thus, we have $\mathbb{P}^n={\rm Proj}(k[X_{1},\ldots,X_{n+1}])$. 
We denote by $\leq$ an arbitrary monomial order on $k[X_{1},\ldots,X_{n+1}]$ such that $X_{1}<\ldots<X_{n+1}$, 
and by $\preceq$ the associated graded lexicographic order, i.e., 
the order in which $X^\gamma\prec X^\delta$ if, and only if, either $|\gamma|<|\delta|$ or 
$|\gamma|=|\delta|$ and $\gamma<\delta$. 
Given a homogeneous ideal $\mathscr{I}$ in $k[X_{1},\ldots,X_{n+1}]$, 
we denote by $\mathscr{C}(\mathscr{I})$ the set of all $\gamma\in\mathbb{N}^{n+1}$ 
which occur as exponents of leading terms of elements of $\mathscr{I}$.
We also consider the complement $\mathscr{D}(\mathscr{I})=\mathbb{N}^{n+1}-\mathscr{C}(\mathscr{I})$, 
which is called the set of standard monomials of $\mathscr{I}$. Again, 
the $k$-vector space $_{k}\langle X^\gamma;\gamma\in \mathscr{D}(\mathscr{I})\rangle$ 
is isomorphic to the $k$-vector space $k[X]/\mathscr{I}$, 
via the canonical map $X^\gamma\mapsto X^\gamma+\mathscr{I}$. 

Let $\mathscr{A}$ be a set of $k$-rational closed points of $\mathbb{P}^n$. 
We may think of the points of $\mathscr{A}$ as lines in $\mathbb{A}^{n+1}$ 
passing through $0$ and at least though one other point of $k^{n+1}$. 
This set defines the ideal $\mathscr{I}(\mathscr{A})$, 
consisting of those polynomials in $k[X_{1},\ldots,X_{n+1}]$ which vanish at all elements of $\mathscr{A}$. 
This is a homogeneous ideal, for which we will consider 
$\mathscr{C}(\mathscr{A})=\mathscr{C}(\mathscr{I}(\mathscr{A}))$ and
$\mathscr{D}(\mathscr{A})=\mathscr{D}(\mathscr{I}(\mathscr{A}))$, in a shorthand notation.
The set $\mathscr{D}(\mathscr{A})$ is not finite, as otherwise, 
ideal $\mathscr{I}$ would only have a finite number of zeros in $L^{n+1}$, for each field extension $L$ of $k$. 
Therefore, $\mathscr{D}(\mathscr{A})$ contains subsets of the form $\gamma+\mathbb{N}e_{i}$, 
for some $\gamma\in\mathbb{N}^{n+1}$ and some $i\in\{1,\ldots,n+1\}$, 
where $e_{i}=(0,\ldots,0,1,0,\ldots,0)$ is the $i$-th standard basis element of $\mathbb{N}^{n+1}$. 
Note that the set $\mathscr{D}(\mathscr{A})\subset\mathbb{N}^{n+1}$ has the property that if 
$\delta$ lies in $\mathscr{D}(\mathscr{A})$, then also $\delta-e_{j}$ lies in $\mathscr{D}(\mathscr{A})$, for all $j$. 
Hence we may assume that the subset $\gamma+\mathbb{N}e_{i}$ of $\mathscr{D}(\mathscr{A})$
is such that the $i$-th component, $\gamma_{i}$, of $\gamma$, is $0$. 
We call a set of the form $\gamma+\mathbb{N}e_{i}$, where $\gamma_{i}=0$, 
an {\it axis} in $\mathbb{N}^{n+1}$. Then the projective analogue of the above affine statement is the following. 
The number of axes in $\mathscr{D}(\mathscr{A})$ equals the number of elements of $\mathscr{A}$.
This will be our main result, Theorem \ref{thma}. 

However, our goal is not only to prove this above statement concerning $\mathscr{D}(\mathscr{A})$, 
but also to construct the reduced Gr\"obner basis $\mathscr{G}(\mathscr{A})$ of the ideal $\mathscr{I}(\mathscr{A})$. 
For reaching both goals, we will proceed from the same basic idea, which is as follows. 
Consider the decomposition of the set $\mathbb{P}^n$ into the two subsets $\mathbb{A}^n\coprod\mathbb{P}^{n-1}$.
Here we think of $\mathbb{P}^n$ in the na\"ive way, i.e. as the set of lines in $k^{n+1}$ passing through $0$. 
In the vector space $k^{n+1}$, we use the coordinates $X_{1},\ldots,X_{n+1}$;
then $\mathbb{A}^n$ is the set of those lines which pass through the subset of $k^{n+1}$ given by $\{X_{1}=1\}$,
a set which can clearly be identified with $k^n$. The ``rest'' of $\mathbb{P}^n$ is then the set of lines 
which lie entirely in the subset of $k^{n+1}$ given by $\{X_{1}=0\}$; this set is identfied with $\mathbb{P}^{n-1}$.
Corresponding to the decomposition $\mathbb{P}^n=\mathbb{A}^n\coprod\mathbb{P}^{n-1}$, 
we divide set $\mathscr{A}$ into the two subsets $\mathscr{A}=\mathscr{A}_{1}\coprod\mathscr{A}_{0}$.
\begin{itemize}
  \item $\mathscr{A}_{1}\subset\mathscr{A}$ consists of lines $\ell\in\mathscr{A}$ such that the intersection 
    $\ell\cap\{X_{1}=1\}$ is nonempty. 
  \item $\mathscr{A}_{0}\subset\mathscr{A}$ consists of lines $\ell\in\mathscr{A}$ such that 
    $\ell\subset\{X_{1}=0\}$. 
\end{itemize}

It follows that intersection $\mathscr{A}_{1}\cap\mathbb{A}^{n}$ is a finite set of closed $k$-rational points of 
$\mathbb{A}^{n}$. We denote this intersection by $A_{1}$. We clearly have $\#A_{1}=\#\mathscr{A}_{1}$. 
Accordingly, let $I(A_{1})$ be the ideal in $k[X_{2},\ldots,X_{n+1}]$
consisting of those polynomials that vanish at all points of $A_{1}$. 
In Section \ref{axes}, we show our main result for $\mathscr{A}_{1}$, 
i.e., we show that $\mathscr{D}(\mathscr{A}_{1})$ contains precisely $\#\mathscr{A}_{1}$ axes. 

In Section \ref{stepone}, we use some of the techniques developed in Section \ref{axes} 
for studying the Gr\"obner basis of $\mathscr{I}(\mathscr{A}_{1})$. 
We present an algorithm for the construction of the graded lexicographic Gr\"obner basis 
$\mathscr{G}_{\rm deglex}(\mathscr{A}_{1})$ of $\mathscr{I}(\mathscr{A}_{1})$, 
starting from the lexicographic Gr\"obner basis $G_{\rm lex}(A_{1})$ of $I(G_{1})$. 
This algorithm in fact works only for the lexicographic (resp. graded lexicographic) order on 
$k[X_{1},\ldots,X_{n}]$, and not for an arbitrary one. 

In Section \ref{allaxes}, we show our main result for general $\mathscr{A}$ by induction over the dimension $n$. 
This is where our basic idea, the decomposition $\mathbb{P}^n=\mathbb{A}^n\coprod\mathbb{P}^{n-1}$, 
and deriving from that, $\mathscr{A}=\mathscr{A}_{1}\coprod\mathscr{A}_{0}$, is being used. 
In Section \ref{axes}, we have shown the statement for $\mathscr{D}(\mathscr{A}_{1})$, and by the induction hypothesis, 
we may assume that the statement is known for $\mathscr{D}(\mathscr{A}_{0})$. 
From these two tokens, we will derive the statement for $\mathscr{D}(\mathscr{A})$. 

In Section \ref{merge}, we will use that same technique for passing from the Gr\"obner bases
$\mathscr{G}(\mathscr{A}_{1})$ and $\mathscr{G}(\mathscr{A}_0)$, both of which we may assume to be known, 
to the Gr\"obner basis $\mathscr{G}(\mathscr{A})$. 

\section{Preliminaries}

Let us first introduce some notation. The leading monomial, resp. the leading exponent 
w.r.t. $\leq$ of $g\in k[X_{1},\ldots,X_{n+1}]$ are denoted by $M(g)$, resp. $E(g)$.
The leading monomial, resp. leading exponent w.r.t. $\preceq$ 
of $h\in k[X_{1},\ldots,X_{n+1}]$ are denoted by $\mathscr{M}(h)$, resp. $\mathscr{E}(h)$.
The total degree of some $\gamma\in\mathbb{N}^{n+1}$, resp. of the corresponding monomial $X^{\gamma}$, 
is denoted by $|\gamma|$, resp. $|X^{\gamma}|$. The total degree of a polynomial is of course the 
maximum of the total degree of its monomials.

In fact, in most cases when we get to work with nonhomogeneous polynomials, 
these will mostly be elements of $I(A_{1})$, thus elements of $k[X_{2},\ldots,X_{n}]$. 
Therefore, we will understand their exponents to lie in $\mathbb{N}^n$, 
viewed as a subset of $\mathbb{N}^{n+1}$ via the embedding
\begin{equation*}
  \mathbb{N}^n\hookrightarrow\mathbb{N}^{n+1}:\alpha\mapsto(0,\alpha)\,.
\end{equation*}
We will frequently work with elements $\alpha$ of $\mathbb{N}^n$.
Whenever we understand $\alpha\in\mathbb{N}^n$ to lie in $\mathbb{N}^{n+1}$ via the above inclusion, 
we denote this element by $(0,\alpha)$. 
(In fact, the only letters we use for elements of $\mathbb{N}^n$ are $\alpha$ and $\beta$, 
and we will not use these letters for elements of $\mathbb{N}^{n+1}$.)

As already mentioned in the Introduction, we will use the sets 
\begin{equation*}
  C(A_{1})=\{E(g);g\in I(A_{1})\}
\end{equation*}
and
\begin{equation*}
  D(A_{1})=\mathbb{N}^{n}-C(A_{1})\,.
\end{equation*}
Furthermore, we denote by $B(A_{1})$ (the ``corners'' of $C(A_{1})$) 
the minimal subset of $C(A_{1})$ such that $C(A_{1})$ is the union of all $\alpha+\mathbb{N}^{n}$, 
where $\alpha$ runs through $B(A_{1})$. In other words, 
$B(A_{1})$ equals the set of lexicographically leading exponents $E(g)$, 
where $g$ runs through the elements of the reduced Gr\"obner basis $G(A_{1})$. Likewise, we set
\begin{equation*}
  \begin{split}
    \mathscr{B}(\mathscr{A}_{1})&=\{\mathscr{E}(h);h\in\mathscr{G}(\mathscr{A}_{1})\}\,,\\
    \mathscr{C}(\mathscr{A}_{1})&=\{\mathscr{E}(h);h\in\mathscr{I}(\mathscr{A}_{1})\}\,,\\
    \mathscr{D}(\mathscr{A}_{1})&=\mathbb{N}^{n+1}-\mathscr{C}(\mathscr{A}_{1})\,,
  \end{split}
\end{equation*}
and make analogous definitions for $\mathscr{A}_{0}$ and $\mathscr{A}_{1}$. 
If there is any necessity to stress the term order $\leq$ (resp. $\preceq$) we are using, 
we add a subscript $_{\leq}$ (resp. $_{\preceq}$) to the objects in question, i.e., 
we write $E_{\leq}$, $M_{\leq}$, $C_{\leq}$, etc. 
(resp. $\mathscr{E}_{\preceq}$, $\mathscr{M}_{\preceq}$, $\mathscr{C}_{\preceq}$ etc.)

We will make use of the projection
\begin{equation*}
  \pi:\mathbb{N}^{n+1}\to\mathbb{N}^{n}:(\gamma_{1},\ldots,\gamma_{n+1})\mapsto(\gamma_{2},\ldots,\gamma_{n+1})\,.
\end{equation*}
The operation corresponding to $\pi$ in the coordinate rings is replacement of 
$h\in k[X_{1},\ldots,X_{n+1}]$ by $h(1,X_{2},\ldots,X_{n+1})\in k[X_{2},\ldots,X_{n+1}]$. 
We denote this polynomial by $h(1,X)$ for short.
We also use an operation in the opposite direction, namely, 
\begin{equation*}
  \mathscr{H}:k[X_{2},\ldots,X_{n+1}]\to k[X_{1},\ldots,X_{n+1}]\,,
\end{equation*}
which maps $g$ to its homogenisation, which is obtained by multiplying each monomial with the smallest power of $X_{1}$
such that the total degree of the product equals $|g|$.

We start with a familiar Lemma; however, for the sake of completeness, 
and due to its great significance to our work, we shall give it a proof. 

\begin{lmm}\label{fund}
  Let $\leq$ be an arbitrary term order in $k[X_{1},\ldots,X_{m}]$, let $I$ be an ideal of $k[X_{1},\ldots,X_{m}]$, 
  and let $C$ the set of leading exponents of elements of $I$, and $D=\mathbb{N}^m-C$. 
  Then for all $\sigma\in C$, there exists a uniquely determined monic polynomial $f_{\sigma}\in I$ such that
  \begin{enumerate}
    \item[(a)] $E(f_{\sigma})$, the leading exponent of $f_{\sigma}$ w.r.t. $\leq$, equals $\sigma$, and
    \item[(b)] all nonleading exponents of $f_{\sigma}$ w.r.t. $\leq$ lie in $D$.
  \end{enumerate}
  Furthermore, the collection of all $f_{\sigma}$, where $\sigma$ runs through $C$, 
  is a basis of the $k$-vector space $I$. 
  In the case where $I$ is a homogeneous ideal, all $f_{\sigma}$ are homogeneous.
\end{lmm}

\begin{proof}
  Existence is shown by induction over $\sigma$. If $\sigma$ is the minimum w.r.t. $\leq$ of $C$, 
  then $f_{\sigma}$ is the unique element of the reduced Gr\"obner basis with leading exponent $\sigma$.
  We assume that the existence of $f_{\sigma^{\prime}}$ is shown for all $\sigma^\prime\in C$ 
  such that $\sigma^{\prime}<\sigma$ and show the existence of $f_{\sigma}$.
  
  If $\sigma$ is the exponent of an element of the reduced Gr\"obner basis of $I$, 
  then we define $f_{\sigma}$ as exactly this element, and $f_{\sigma}$ satisfies (a) and (b).
  If $\sigma$ is not the exponent of an element of the reduced Gr\"obner basis, 
  then there exists an index $i\in\{1,\ldots,m\}$ such that $\sigma-e_{i}$ lies in $C$. 
  We take such an $i$ and define $\sigma^{\prime}=\sigma-e_{i}$. 
  Clearly $\sigma^{\prime}<\sigma$, hence there exists $f_{\sigma^{\prime}}$ as claimed, 
  and $E(X_{i}f_{\sigma^\prime})=\sigma$. If all nonleading exponents of $X_{i}f_{\sigma^\prime}$ lie in $D$, 
  we take $f_{\sigma}=X_{i}f_{\sigma^\prime}$. Otherwise, let $\sigma^{\prime\prime}$ 
  be the largest nonleading exponent of the polynomial $f_{\sigma}^{(0)}=X_{i}f_{\sigma^\prime}$ 
  which lies in $C$, and let $c$ be the coefficient of $X^{\sigma^{\prime\prime}}$ in $f_{\sigma}^{(0)}$. 
  We replace $f_{\sigma}^{(0)}$ by $f_{\sigma}^{(1)}=f_{\sigma}^{(0)}-cf_{\sigma^{\prime\prime}}$. 
  This polynomial also lies in $I$, and $E(f_{\sigma}^{(1)})=\sigma$, 
  but the smallest nonleading term of $f_{\sigma}^{(1)}$ which lies in $C$ (if it exists)
  is smaller than the smallest nonleading term of $f_{\sigma}^{(0)}$ which lies in $C$. 
  We repeat the last argument, thus constructing $f_{\sigma}^{(2)},f_{\sigma}^{(3)},\ldots$, 
  then this process will at some point terminate, and we will finally find a polynomial satisfying (a) and (b).
  
  As for the uniqueness of $f_{\sigma}$, let us assume that also $f_{\sigma}^{\prime}$ satisfies (a) and (b). 
  Then the leading exponent of $f_{\sigma}-f_{\sigma}^{\prime}$ lies in $D$, 
  hence $f_{\sigma}-f_{\sigma}^{\prime}=0$.
  
  Linear independence of $f_{\sigma}$, where $\sigma$ runs through $C$, 
  follows from linear independence of the respective leading monomials $X^\sigma$. 
  Let us show that the span of $f_{\sigma}$, where $\sigma$ runs through $C$, is the whole $k$-vector space $I$. 
  Given $f\in C$, consider $\sigma^{(1)}$, the leading exponent of $f$, and $c^{(1)}$ be the leading coefficient of $f$. 
  Then the polynomial $f^{(1)}=f-c^{(1)}f_{\sigma^{(1)}}$ also lies in $I$, but $E(f^{(1)})<E(f)$. 
  Analogously, we construct $f^{(2)},f^{(3)},\ldots$; this process will terminate, 
  and $f$ indeed lies in the span of all $f_{\sigma}$. 
  
   If $I$ is a homogeneous ideal, then along with each element $f$ of $I$, 
   also all homogeneous components of $f$ lie in $I$. 
  Let $f_{\sigma}$ satisfy (a) and (b), and denote by $f_{\sigma}=f_{\sigma,0}+\ldots+f_{\sigma,d}$ 
  the decomposition of $f_{\sigma}$ into homogeneous elements of degree $0,\ldots,d$. 
  For $j<d$, we have $f_{\sigma,j}\in I$, and $E(f_{\sigma})\in D$. 
  Therefore, $f_{\sigma,j}=0$ for all $j<d$, and $f_{\sigma}$ is indeed homogeneous.
\end{proof}

\section{Axes in $\mathscr{D}(\mathscr{A}_{1})$}\label{axes}

Let us now fix an arbitrary term order $\leq$ on the polynomial ring $k[X_{1},\ldots,X_{n+1}]$, 
such that $X_{1}<\ldots<X_{n}$. 
This term order clearly induces a term order on the subring $k[X_{2},\ldots,X_{n+1}]$. 

\begin{lmm}\label{span} 
  The ideal $\mathscr{I}(\mathscr{A}_{1})$ is the $k$-span of 
  \begin{equation*}
    X_{1}^{\epsilon}\mathscr{H}(f)\,,
  \end{equation*}    
    where $\epsilon$ runs through $\mathbb{N}$, and $f$ runs through $I(A_{1})$. 
    Furthermore, each $h\in\mathscr{G}(\mathscr{A}_{1})$ is the homogenisation of some $g\in I(A_{1})$. 
\end{lmm}

\begin{proof}
  For the first assertion, let $h\in\mathscr{I}(\mathscr{A}_{1})$ be given. 
  We may assume that $h$ is homogeneous, as all homogeneous components of a given element of 
  $\mathscr{I}(\mathscr{A}_{1})$ lie in $\mathscr{I}(\mathscr{A}_{1})$.
  Let $X_{1}^{\epsilon}$ be the smallest power of $X_{1}$ which divides $h$. Write 
  \begin{equation*}
    h=X_{1}^{\epsilon}h^\prime\,. 
  \end{equation*}
  Definition of $\epsilon$ implies that $h^\prime=\mathscr{H}(h^\prime(1,X))$. 
  On the other hand, $h^\prime(1,X)\in I(A_{1})$, as $h^\prime(1,X)=h(1,X)$. 
  
  For the second assertion, let $h$ be an element of $\mathscr{G}(\mathscr{A}_{1})$.
  We decompose $h=h_{0}+\ldots+h_{d}$ into its homogeneous components, 
  all of which lie in $\mathscr{I}(\mathscr{A}_{1})$. 
  Consider a summmand $h_{j}$, where $j<d$. We first show that $h$ is homogeneous, i.e., $h_{j}=0$.
  Assume that $h_{j}\neq0$. By the very definition of the Gr\"obner basis, 
  there exists an element $h^\prime$ of $\mathscr{G}(\mathscr{A}_{1})$
  such that $\mathscr{M}(h^{\prime})$ divides $\mathscr{M}(h_{j})$. 
  This particular $h^{\prime}$ cannot equal $h$, as
  \begin{equation*}
    |\mathscr{M}(h^{\prime})|\leq|\mathscr{M}(h_{j})|<|\mathscr{M}(h_{d})|=|\mathscr{M}(h)|\,.
  \end{equation*}
  Furthermore, the fact that the Gr\"obner basis $\mathscr{G}(\mathscr{A})$ is {\it reduced} implies that
  no monomial occurring in $h$ is divided by any $\mathscr{M}(h^{\prime})$, 
  where $h^{\prime}$ runs through $\mathscr{G}(\mathscr{A})-\{g\}$. Thus in particular, for our $h_{j}$, 
  we get that $\mathscr{M}(h^{\prime})\nmid\mathscr{M}(h_{j})$, for all $h^{\prime}\in\mathscr{G}(\mathscr{A})-\{g\}$.
  Hence $\mathscr{M}(h^{\prime})\nmid\mathscr{M}(h_{j})$, for all $h^{\prime}\in\mathscr{G}(\mathscr{A})$.
  This is a contradiction, which yields $h_{j}=0$, as claimed.
  
  Let $X_{1}^{\epsilon}$ be the largest power of $X_{1}$ which divides $h$, 
  and $h^{\prime}$ be the other factor in $h$, i.e., $h=X_{1}^{\epsilon}h^{\prime}$. Then as $h(1,X)=h^{\prime}(1,X)$, 
  polynomial $h^{\prime}$ also lies in $\mathscr{I}(\mathscr{A}_{1})$. 
  As $\mathscr{E}(h)=\mathscr{E}(h^{\prime})+\epsilon e_{1}$, and as $h\in\mathscr{G}(\mathscr{A}_{1})$, 
  we have $\epsilon=0$. From this follows the second assertion.
\end{proof}

Next, we wish to gain some information on $\pi(\mathscr{C}(\mathscr{A}_{1}))$.
However, this will be feasible only in that particular case where the order $\leq$ is the lexicographic order, 
hence $\preceq$ is the graded lexicographic order. 
We denote these orders by $\leq_{\rm lex}$, resp. $\preceq_{\rm deglex}$. 
Note that $X_{1}<_{\rm lex}\ldots<_{\rm lex}X_{n+1}$, and $X_{1}\prec_{\rm deglex}\ldots\prec_{\rm deglex}X_{n+1}$.

\begin{lmm}\label{lead}
  Let $h\in k[X_{1},\ldots,X_{n+1}]$ be a homogeneous polynomial, then
  \begin{equation*}
    E_{\rm lex}(h(1,X))=\mathscr{E}_{\rm deglex}(h(1,X))=\pi(\mathscr{E}_{\rm deglex}(h))\,.
  \end{equation*}
\end{lmm}

\begin{proof}
  Let $\gamma$ and $\delta$ be exponents of monomials of $h$. Then $|\gamma|=|\delta|$. 
  This implies that $\gamma\prec_{\rm deglex}\delta$ if, and only if, $\gamma<_{\rm lex}\delta$. 
  The latter condition clearly implies that $\pi(\gamma)\leq_{\rm lex}\pi(\delta)$. 
  In this equality, the case $\pi(\gamma)=\pi(\delta)$ is impossible, 
  since from it, together with $|\gamma|=|\delta|$, we would get $\gamma=\delta$. 
  Therefore, condition $\gamma<_{\rm lex}\delta$
  is also equivalent to $\pi(\gamma)<_{\rm lex}\pi(\delta)$. This proves the Lemma.
\end{proof}

Note that we have used a specific property of the lexicographic order 
when deducing from $\gamma<_{\rm lex}\delta$ that $\pi(\gamma)\leq_{\rm lex}\pi(\delta)$. 
Here is an example of a term order for which this fails. 
Take $\leq$ to be the graded lexicographic order 
(hence the associated $\preceq$ is also the graded lexicographic order), 
$\gamma=(0,2)$ and $\delta=(2,1)$. 

\begin{pro}\label{c}
  $\pi(\mathscr{C}_{\rm deglex}(\mathscr{A}_{1}))=C_{\rm lex}(A_{1})$.
\end{pro}

\begin{proof}
  Given $\alpha\in C_{\rm lex}(A_{1})$, polynomial $f_{\alpha}$ gives rise to $\mathscr{H}(f_{\alpha})$, 
  which lies in $\mathscr{I}(\mathscr{A}_{1})$, 
  and whose leading exponent w.r.t. $\prec$ is $(\epsilon,\alpha)$, 
  for a suitable $\epsilon\in\mathbb{N}$. 
  This proves inclusion $\pi(\mathscr{C}_{\rm lex}(\mathscr{A}_{1}))\supset C_{\rm lex}(A_{1})$.
  
  For inclusion $\pi(\mathscr{C}_{\rm deglex}(\mathscr{A}_{1}))\subset C_{\rm lex}(A_{1})$, 
  we take $h\in\mathscr{I}(\mathscr{A}_{1})$ and show that
  $\pi(\mathscr{E}_{\rm deglex}(h))$ lies in $C_{\rm lex}(A_{1})$. 
  We may assume that $h$ is homogeneous. 
  By Lemma \ref{lead}, we have $\pi(\mathscr{E}_{\rm deglex}(h))=E_{\rm lex}(h(1,X))$; 
  as $h(1,X)\in I(A_{1})$, exponent $E_{\rm lex}(h(1,X))$ lies in $C_{\rm lex}(A_{1})$.
\end{proof}

\begin{cor}\label{corry}
  $\pi(\mathscr{B}_{\rm deglex}(\mathscr{A}_{1}))\subset C_{\rm lex}(A_{1})$.
\end{cor}

Let us now consider axes in $\mathbb{N}^{n+1}$, as defined in the Introduction, i.e., sets of the form
$\gamma+\mathbb{N}e_{i}$ for some fixed $\gamma\in\mathbb{N}^{n+1}$
such that $\gamma_{i}=0$. An axis of this form is called {\it parallel to $\mathbb{N}e_{i}$}.
We understand an axis to be a subobject of $\mathbb{N}^{n+1}$ of dimension one;
in contrast to that, an element $\gamma\in\mathbb{N}^{n+1}$ is understood to be a subobject of dimension zero. 
Proposition \ref{c} enables us to determine what 
$\mathscr{D}_{\rm deglex}(\mathscr{A}_{1})$ looks like from ``bird's eye view'', 
i.e. to determine all axes contained in $\mathscr{D}_{\rm deglex}(\mathscr{A}_{1})$. 
The zero-dimensional subobjects of $\mathscr{D}_{\rm deglex}(\mathscr{A}_{1})$ will be computed later. 

\begin{thm}\label{thmone}
  $\mathscr{D}_{\rm deglex}(\mathscr{A}_{1})$ contains precisely $\#D_{\rm lex}(A_{1})$ axes. These are given by
  \begin{equation*}
    (0,\alpha)+\mathbb{N}e_{1}\,,\text{ for }\alpha\in D_{\rm lex}(A_{1})\,.
  \end{equation*}
\end{thm}

\begin{proof}
  The existence of the axes as claimed follows from the identity 
  $\pi(\mathscr{C}_{\rm deglex}(\mathscr{A}_{1}))=C_{\rm lex}(A_{1})$
  of Proposition \ref{c}. Likewise, this identity implies that $\mathscr{C}_{\rm deglex}(\mathscr{A}_{1})$ 
  contains no other axes parallel to $\mathbb{N}e_{1}$ than those claimed.
  It remains to show that $\mathscr{D}_{\rm deglex}(\mathscr{A}_{1})$ contains no axes parallel to $\mathbb{N}e_{j}$, 
  for all $j>1$. For this, it suffices to show that, given $j>1$, 
  there exists $h\in\mathscr{I}(\mathscr{A}_{1})$ such that $\mathscr{E}_{\rm deglex}(h)$ lies in $\mathbb{N}e_{j}$. 
  For $j$ as above, take $\mu\in\mathbb{N}$ strictly larger than $|\gamma|$, 
  for all $\gamma\in D_{\rm lex}(A)$. Set $\alpha=\mu e_{j}$; it follows that $|f_{\alpha}|=|\alpha|$, 
  as all lexicographically nonleading exponents of $f_{\alpha}$ lie in $D_{\rm lex}(A)$. Therefore, 
  the graded lexicographically leading exponent of $h=\mathscr{H}(f_{\alpha})$ is $\alpha$, 
  which indeed lies on the axis $\mathbb{N}e_{j}$. 
\end{proof}

The first assertion of Theorem \ref{thmone} holds for an arbitrary term order $\leq$. 
For proving this, we need to prove two further facts.

\begin{lmm}\label{direction}
  All axes contained in $\mathscr{D}(\mathscr{A}_{1})$ are parallel to $\mathbb{N}e_{1}$. 
\end{lmm}

\begin{proof}
  We show that on each axis $\mathbb{N}e_{j}$, where $j>1$, 
  lies the leading exponent of an element of $\mathscr{I}(\mathscr{A}_{1})$. 
  The assertion follows from this. Denote by $s$ the cardinality of $\mathscr{A}_{1}$;
  thus we can think of $\mathscr{A}_{1}$ as the union of lines $\ell_{i}$ in $\mathbb{A}^{n+1}$, 
  where $i=1,\ldots,s$; as the elements of $\mathscr{A}_{1}$ are $k$-rational, 
  line $\ell_{i}$ is given by a system of equations
  \begin{equation*}
    X_{j}-a_{i,j}X_{1}=0\,,\text{ for }j\in\{2,\ldots,n+1\}\,.
  \end{equation*}
  Fix $j>1$ and consider the element
  \begin{equation*}
    p=\prod_{i=1}^s(X_{j}-a_{i,j}X_{1})
  \end{equation*}
  of $\mathscr{I}(\mathscr{A}_{1})$. When muliplying all factors of $p$, 
  we find that $p$ is a sum of $X_{j}^s$ plus $k$-multiples of monomials $X_{1}^{\gamma_{1}}X_{j}^{\gamma_{j}}$, 
  where $\gamma_{1}+\gamma_{j}=s$ and $\gamma_{j}<s$. 
  The latter monomials are all strictly smaller than $X_{j}^s$, as
  \begin{equation*}
    X_{1}^s=X_{1}^{s-1}X_{1}<X_{1}^{s-1}X_{j}=X_{1}^{s-2}X_{1}X_{j}<X_{1}^{s-2}X_{j}^2<\ldots<X_{j}^s\,.
  \end{equation*}
  Thus $\mathscr{E}(p)=se_{j}$, as claimed. 
\end{proof}

The following Proposition holds without the assumption that the term order $\preceq$ is the graded order
associated to $\leq$; therefore, we formulate the Proposition in this generality.

\begin{pro}\label{parallel}
  Let $\preceq$ and $\sqsubseteq$ be two term orders on $k[X_{1},\ldots,X_{n+1}]$
  such that $X_{1}\prec\ldots\prec X_{n+1}$ and $X_{1}\sqsubset\ldots\sqsubset X_{n+1}$. 
  Let $I$ be an ideal in $k[X_{1},\ldots,X_{n+1}]$ such that all axes contained in 
  $\mathscr{D}_{\preceq}(I)$ and all axes contained in $\mathscr{D}_{\sqsubseteq}(I)$
  are parallel to $\mathbb{N}e_{1}$. Denote by $N_{\preceq}(I)$ (resp. $N_{\sqsubseteq}(I)$) the number of 
  axes contained in $\mathscr{D}_{\preceq}(I)$ (resp. $\mathscr{D}_{\sqsubseteq}(I)$). 
  Then 
  \begin{equation*}
    N_{\preceq}(I)=N_{\sqsubseteq}(I)\,.
  \end{equation*}
\end{pro}

\begin{proof}
  We shall prove that $N_{\preceq}(I)\leq N_{\sqsubseteq}(I)$. 
  As the roles of term orders $\preceq$ and $\sqsubseteq$ are symmetric, this will prove the Proposition.

  Consider the decomposition 
  \begin{equation*}
    D_{\preceq}(I)=Y_{\preceq}\coprod Y_{\preceq}^\prime\,,
  \end{equation*}
  where $Y_{\preceq}=\pi(Y_{\preceq})\times\mathbb{N}e_{1}$, and $Y_{\preceq}^\prime$ is a finite set. 
  Thus $Y_{\preceq}$ consists of all affine axes contained in $D_{\preceq}(I)$, 
  and $Y_{\preceq}^\prime$ is the complement of $Y_{\preceq}$ in $\mathscr{D}_{\preceq}(I)$.
  Analogously, we define the decomposition 
  \begin{equation*}
    D_{\sqsubseteq}(I)=Y_{\sqsubseteq}\coprod Y_{\sqsubseteq}^\prime\,.
  \end{equation*}
  The $k$-vector space $R_{\preceq}=\oplus_{\gamma\in D_{\preceq}(I)}k X^\gamma$ 
  is a system of representatives of the $k$-vector space $k[X]/I$.
  We define a map $r_{\preceq}:k[X]\to R_{\leq}$ by sending a polynomial $f$ to the unique representative of 
  $f+I$ in $R_{\leq}$. 
  This induces an isomorphism of $k$-vector spaces, which we denote by the same letter, 
  $r_{\preceq}:k[X]/I\to R_{\preceq}$.
  Analogously, we define a set $R_{\sqsubseteq}$, a map $r_{\sqsubseteq}:k[X]\to R_{\sqsubseteq}$, 
  and an isomorphism of $k$-vector spaces, $r_{\sqsubseteq}:k[X]/I\to R_{\sqsubseteq}$. 
  
  Now we consider a subset $M$ of $Y_{\sqsubseteq}$ of the form 
  $M=\pi(Y_{\sqsubseteq})\times\{0,\ldots,\mu\}e_{1}$, for some $\mu\geq0$, 
  such that for all $\alpha\in\pi(Y_{\preceq})$
  and for all $\gamma\in (e_{1}+Y_{\sqsubseteq}^\prime)$, we have 
  \begin{equation}\label{l}
    r_{\sqsubseteq}(X^\alpha)\,,r_{\sqsubseteq}(X^\gamma)
    \in_{k}\langle X^\delta;\delta\in M\coprod Y_{\sqsubseteq}^\prime\rangle\,.
  \end{equation}
  (As this condition on $M$ involves only a finite number $\alpha$ and $\gamma$, 
  a set $M$ of the form $M=\pi(Y_{\sqsubseteq})\times\{0,\ldots,\mu\}e_{1}$ indeed exists.) 
  Let us now write an element $\gamma$ of $Y_{\preceq}$ as $\gamma=(\epsilon,\pi(\gamma))$. 
  We claim that 
  \begin{equation}\label{le}
    r_{\sqsubseteq}(X^\gamma)\in_{k}\langle X^\delta;\delta\in M_{\epsilon}\coprod Y_{\sqsubseteq}^\prime\rangle\,,
  \end{equation}
  where 
  \begin{equation*}
    M_{\epsilon}=\cup_{0\leq\epsilon^\prime\leq\epsilon}(M+\epsilon^\prime e_{1})\,.
  \end{equation*}
  We show this by induction over $\epsilon$. The case $\epsilon=0$ is clear by definition of $M=M_{0}$. 
  As for the induction step, let $\gamma=(\epsilon,\pi(\gamma))$ be given. 
  We may assume that the claim is true for $\gamma^\prime=\gamma-e_{1}=(\epsilon-1,\pi(\gamma))$. 
  Therefore, there exist coefficients $c_{\delta}\in k$ such that
  \begin{equation*}
    r_{\sqsubseteq}(X^{\gamma^\prime})
    =\sum_{\delta\in M_{\epsilon-1}}c_{\delta}X^\delta+\sum_{\delta\in Y_{\sqsubseteq}^\prime}c_{\delta}X^\delta\,.
  \end{equation*}
  It follows that $r_{\sqsubseteq}(X^{\gamma})$ takes the shape
  \begin{equation*}
    r_{\sqsubseteq}(X^{\gamma})=r_{\sqsubseteq}(X_{1}X^{\gamma^\prime})
    =r_{\sqsubseteq}(\sum_{\delta\in M_{\epsilon-1}}c_{\delta}X_{1}X^\delta
    +\sum_{\delta\in Y_{\sqsubseteq}^\prime}c_{\delta}X_{1}X^\delta)\,.
  \end{equation*}
  The inclusion $M_{\epsilon}\subset Y_{\sqsubseteq}$ implies that if $\delta$ lies in $M_{\epsilon-1}$, 
  the element $\delta+e_{1}$ lies in $Y_{\sqsubseteq}$. Hence, 
  $r_{\sqsubseteq}(\sum_{\delta\in M_{\epsilon-1}}c_{\delta}X_{1}X^\delta)
  =\sum_{\delta\in M_{\epsilon-1}}c_{\delta}X_{1}X^\delta$, 
  which is contained in the $k$-vector space on the right hand side of \eqref{le}. 
  Furthermore, by definition of $M$, 
  sum $\sum_{\delta\in Y_{\sqsubseteq}^\prime}c_{\delta}X_{1}X^\delta$ has a representant in
  the $k$-vector space on the right hand side of \eqref{l}. As $M\subset M_{\epsilon}$, 
  this implies that $r_{\sqsubseteq}(\sum_{\delta\in Y_{\sqsubseteq}^\prime}c_{\delta}X_{1}X^\delta)$ 
  is contained in the $k$-vector space on the right hand side of \eqref{le}. The claim is proved.

  Let $\gamma$ be an arbitrary element of $\pi(Y_{\preceq})\times\{0,\ldots,\epsilon\}e_{1}$. 
  Then $r_{\preceq}(X^\gamma)=X^\gamma$, 
  hence the image of $X^\gamma$ under the composition of $k$-vector space isomorphisms, 
  \begin{equation*}
    \begin{CD}
      r_{\sqsubseteq}\circ r_{\preceq}^{-1}:R_{\preceq}@>r_{\preceq}^{-1}>>k[X]/I@>r_{\sqsubseteq}>> R_{\sqsubseteq}\,,
    \end{CD}
  \end{equation*}
  is $r_{\sqsubseteq}(X^\gamma)$. 
  When letting $\gamma$ run through the set $\pi(Y_{\preceq})\times\{0,\ldots,\epsilon\}e_{1}$, we get a collection of 
  \begin{equation*}
    \#\pi(Y_{\preceq})(\epsilon+1)=N_{\preceq}(I)(\epsilon+1) 
  \end{equation*}
  linearly independent elements of $R_{\preceq}$. Therefore, the collection $r_{\sqsubseteq}(X^\gamma)$, 
  where $\gamma$ runs through $\pi(Y_{\preceq})\times\{0,\ldots,\epsilon\}e_{1}$, 
  is a collection of linearly independent elements of $R_{\sqsubseteq}$. 
  By \eqref{le}, all these elements are contained in the subset 
  $_{k}\langle X^\delta;\delta\in Y_{\sqsubseteq}^\prime\coprod M_{\epsilon}\rangle$ of $R_{\sqsubseteq}$. 
  This set contains 
  \begin{equation*}
    \#Y_{\sqsubseteq}^\prime+\#M_{\epsilon}=\#Y_{\sqsubseteq}^\prime+\#\pi(Y_{\sqsubseteq})(\epsilon+\mu+1)
    =\#Y_{\sqsubseteq}^\prime+N_{\sqsubseteq}(I)(\epsilon+\mu+1)
  \end{equation*}
  linearly independent elements. Therefore, we get
  \begin{equation*}
    N_{\preceq}(I)(\epsilon+1)\leq\#Y_{\sqsubseteq}^\prime+N_{\sqsubseteq}(I)(\epsilon+\mu+1)
  \end{equation*}
  for all $\epsilon\geq0$. It follows that $N_{\preceq}(I)\leq N_{\sqsubseteq}(I)$. 
\end{proof}

\begin{thm}\label{thma1}
  $\mathscr{D}(\mathscr{A}_{1})$ contains precisely $\#D(A_{1})$ axes.
\end{thm}

\begin{proof}
  This follows from Theorem \ref{thmone}, Lemma \ref{direction} and Proposition \ref{parallel}, 
  (where $\sqsubseteq$ is the graded lexicographic order $\preceq_{\rm deglex}$), 
  and the fact that the number $\#D(A_{1})$ is independent of the term order, as it equals $\#A_{1}$.
\end{proof}

\section{The ideal of $\mathscr{A}_{1}$ for the lexicographic order}\label{stepone}

Now we have a good impression of the shape of $\mathscr{D}(\mathscr{A}_{1})$,
knowing all its one-dimensional subobjects. 
Unfortunately, we cannot give an a priori description of the zero-dimensional subobjects of $\mathscr{D}(\mathscr{A}_{1})$
in a state of explicitness comparable to Theorem \ref{thma1}, 
or to Theorem \ref{thmone} in the special case where $\preceq$ is the graded lexicographic order.
However, there exists a series of articles dedicated to the construction of algorithms for the Gr\"obner basis 
$\mathscr{G}(\mathscr{A}_{1})$, see \cite{mora}, \cite{mmm}, \cite{akr1} and \cite{akr2}. 
These algorithms yield, in particular, an algorithm for the determination of $\mathscr{D}(\mathscr{A}_{1})$.
In this sense, we do have knowlegde of those elements of $\mathscr{D}(\mathscr{A}_{1})$ 
that do not occurr in any axis conained in $\mathscr{D}(\mathscr{A}_{1})$. 

In this Section we will give yet another algorithm for the determination of Gr\"obner basis $\mathscr{G}(\mathscr{A}_{1})$, 
which in fact works solely for the special case where the term orders $\leq$ and $\preceq$ are
$\leq_{\rm lex}$ and $\preceq_{\rm deglex}$.
So in fact, we only give an algorithm for $\mathscr{G}_{\rm deglex}(\mathscr{A}_{1})$. 
This algorithm will be based on the ideas developed in Section \ref{axes}, 
more specifically all ingredients of Theorem \ref{thmone}. 
However, the passage from $\mathscr{G}_{\rm deglex}(\mathscr{A}_{1})$ to the Gr\"obner basis
w.r.t. an arbitrary term order is subject to the theories of universal Gr\"obner basis, and of the Gr\"obner fan, 
which are well-known, see \cite{mora}, \cite{weis}, \cite{sturmfels}, and \cite{using}. 
Therefore, we shall not pursue the development of an algorithmic analogue of the proof of Theorem \ref{thma1}.

As our algorithm for the Gr\"obner basis works only for the lexicographic order, 
the present Section might be considered merely a special case of the articles cited above, therefore obsolete. 
However, it is of interest to see how the results of Section \ref{axes} behave in an algorithmic setting.

Thoughout the construction of $\mathscr{G}_{\rm deglex}(\mathscr{A}_{1})$, we will make frequent use of the $k$-basis
$f_{\alpha}$, where $\alpha$ runs through $C_{\rm lex}(A_{1})$, of $I(A_{1})$. 
It is indeed justified to assume that we have this basis at hand.
Firstly, the Gr\"obner basis of ideal $I(A_{1})$ has been studied by many authors
and is sufficiently well-known, also in term of its algorithmic construction. 
In fact, this Gr\"obner basis is subject of the Buchberger--M\"oller algorithm, see the original article \cite{buchmoell}, 
or \cite{big}, or else \cite{jpaa} and the references therein, for the special case where solely the lexicographic order is used. 
Secondly, the proof of Lemma \ref{fund} gives an algorithm for the construction of all
 $f_{\alpha}$, where $\alpha\in C_{\rm lex}(A_{1})$, when given the reduced Gr\"obner basis of $I(A_{1})$. 

By Lemma \ref{span}, 
each $h\in\mathscr{G}_{\rm deglex}(\mathscr{A}_{1})$ arises as homogenisation of some $g\in I(A_{1})$, 
which takes the form 
\begin{equation}\label{g}
  g=f_{\alpha}+\sum_{\substack{\beta\in C_{\rm lex}(A_{1})\\\beta<_{\rm lex}\alpha}}c_{\beta}f_{\beta}\,.
\end{equation}
The leading exponent w.r.t. $\preceq_{\rm deglex}$ of $h$ is $(\epsilon,\alpha)$, where 
\begin{equation*}
  \epsilon=|g|-|\alpha|\,. 
\end{equation*}
(At this point, we have used the following fact about the lexicographic order. 
If $\alpha$ and $\alpha^\prime$ are in $\mathbb{N}^n$ such that $X^\alpha<_{\rm lex}X^{\alpha^\prime}$, 
then also $X_{1}^\epsilon X^\alpha<_{\rm lex}X_{1}^\mu X^{\alpha^\prime}$ for all natural numbers $\epsilon$ and $\mu$.)
When replacing $g$ by $h=\mathscr{H}(g)$, 
each monomial $X^{\alpha^\prime}$ which appears in $g$ is replaced by
$X_{1}^{\epsilon(\alpha^\prime)}X^{\alpha^\prime}$, where $\epsilon(\alpha^\prime)=|g|-|\alpha^\prime|$. 
In particular, we have $\epsilon=\epsilon(\alpha)$ in this notation.

Here are more conditions about the shape a polynomial $g\in I(A_{1})$ has to have so that 
$h=\mathscr{H}(g)$ lies in $\mathscr{G}(\mathscr{A}_{1})$. 

\begin{pro}\label{onetwo}
  Let $\alpha\in C_{\rm lex}(A_{1})$ be given. 
  The homogenisation $h=\mathscr{H}(g)$ of a polynomial $g\in I(A_{1})$ 
  of the form \eqref{g} appears in $\mathscr{G}_{\rm deglex}(\mathscr{A}_{1})$ if, and only if, 
  \begin{enumerate}
    \item[{\it (i)}] $|g|$ is minimal amongst the total degrees of all polynomials of the form \eqref{g}, 
    \item[{\it (ii)}] for all $\beta<_{\rm lex}\alpha$ in $C_{\rm lex}(A_{1})$ and for all 
      $\alpha^{\prime}\in\pi(\mathscr{B}_{\rm deglex}(\mathscr{A}_{1}))$, 
      $\alpha^{\prime}\neq\alpha$, such that $\beta\in\alpha^{\prime}+\mathbb{N}^{n}$, we have 
      \begin{equation*}
        |g|-|\beta|<|g^{\prime}|-|\alpha^{\prime}|\,,
      \end{equation*}
      unless $c_{\beta}=0$. Here, $|g^{\prime}|$ is the minimum of the total degrees of all polynomials of the form 
      \begin{equation}\label{gprime}
        g^{\prime}=f_{\alpha^{\prime}}
        +\sum_{\substack{\beta^\prime\in C_{\rm lex}(A_{1})\\\beta^{\prime}<_{\rm lex}\alpha^{\prime}}}
        c^{\prime}_{\beta^{\prime}}f_{\beta^{\prime}}\,.
      \end{equation}
  \end{enumerate}
  A polynomial $g$ of the form \eqref{g} satisfying (i) and (ii), if it exists, is uniquely determined by these conditions. 
\end{pro}

\begin{proof}
  Let us first show necessity of (i). Assume that there exists some 
  $g^{\prime\prime}=f_{\alpha}+\sum_{\beta<_{\rm lex}\alpha}c_{\beta}^{\prime\prime}f_{\beta}$ 
  such that $|g^{\prime\prime}|<|g|$ and 
  $h^{\prime\prime}=\mathscr{H}(g^{\prime\prime})\in\mathscr{G}_{\rm deglex}(\mathscr{A}_{1})$. 
  Then the leading exponent w.r.t. $\preceq_{\rm deglex}$ of $h^{\prime\prime}$ 
  is $(\epsilon^{\prime\prime},\alpha)$, where $\epsilon^{\prime\prime}<\epsilon$. 
  This is impossible, since we assume $h$ to lie in a Gr\"obner basis which is reduced, hence, in particular, minimal. 
  Necessity of (i) is shown. Next, let us show necessity of (ii). 
  
  The Gr\"obner basis $\mathscr{G}_{\rm deglex}(\mathscr{A}_{1})$ is reduced. 
  This means that each $h\in\mathscr{G}_{\rm deglex}(\mathscr{A}_{1})$ is monic,
  and no term of $h$ is divided by the leading term of any $h^{\prime}\in\mathscr{G}_{\rm deglex}(\mathscr{A}_{1})$, 
  $h^{\prime}\neq h$. Given $g$ as in \eqref{g}, we have seen that its leading term is $X_{1}^{\epsilon}X^{\alpha}$, 
  where $\epsilon=\epsilon(\alpha)=|g|-|\alpha|$.
  There are two classes of nonleading terms occurring in $h$.
  \begin{itemize}
    \item For $\beta<_{\rm lex}\alpha$, if $c_{\beta}\neq0$, 
      then $h$ contains the summand $c_{\beta}X_{1}^{\epsilon(\beta)}X^{\beta}$. 
      (Whenever $c_{\beta}\neq0$, summand $c_{\beta}X^{\beta}$ appears in $g$.
      This is due to the linear independence of $X^\beta$, the leading term of $f_{\beta}$, 
      where $\beta$ runs through $C_{\rm lex}(A_{1})$.)
    \item All other nonleading terms of $h$ are $k$-multiples of homogenisations of $X^{\gamma}$, 
      where $\gamma\in D_{\rm lex}(A_{1})$.
  \end{itemize}
  Analogously as above, if $h^{\prime}\neq h$ is another element of $\mathscr{G}_{\rm deglex}(\mathscr{A}_{1})$, 
  it is the homogenisation of a polynomial $g^\prime$ of the form \eqref{gprime}, 
  where $\alpha^{\prime}\in\pi(\mathscr{B}_{\rm deglex}(\mathscr{A}_{1}))$, 
  subject to the minimality condition on its total degree. 
  The leading term of $h^{\prime}$ is 
  $X_{1}^{\epsilon^{\prime}}X^{\alpha^{\prime}}$, where $\epsilon^{\prime}=|g^{\prime}|-|\gamma^{\prime}|$.
  We have to make sure that $X_{1}^{\epsilon^{\prime}}X^{\alpha^{\prime}}$ divides neither 
  $X_{1}^{\epsilon(\beta)}X^{\beta}$, for $\beta\in C_{\rm lex}(A_{1})$ 
  such that $\beta<_{\rm lex}\alpha$ and $c_{\beta}\neq0$, nor those monomials of $h$ 
  appearing as homogenisation of $X^{\gamma}$, where $\gamma\in D_{\rm lex}(A_{1})$. However, 
  if $X_{1}^{\epsilon^{\prime}}X^{\alpha^{\prime}}$ were to divide an $X_{1}$-multiple of $X^{\gamma}$, 
  this would in particular imply that 
  $\gamma\in\pi((\epsilon^{\prime},\alpha^{\prime}))+\mathbb{N}^{n}=\alpha^{\prime}+\mathbb{N}^{n}$. 
  But this is impossible, as $\gamma\in D_{\rm lex}(A)$, 
  and $\alpha^{\prime}\in C_{\rm lex}(A_{1})$ by Corollary \ref{corry}. 
  
  Hence we only have to make sure that $X_{1}^{\epsilon^{\prime}}X^{\alpha^{\prime}}$ does not divide 
  $X_{1}^{\epsilon(\beta)}X^{\beta}$, for $\beta\in C_{\rm lex}(A_{1})$ 
  such that $\beta<_{\rm lex}\alpha$ and $c_{\beta}\neq0$. 
  In terms of exponents, we have to make sure that 
  $(\epsilon(\beta),\beta)\notin(\epsilon^{\prime},\alpha^{\prime})+\mathbb{N}^{n+1}$ 
  for all $\beta<_{\rm lex}\alpha$ such that $c_{\beta}\neq0$. 
  If inclusion $(\epsilon(\beta),\beta)\in(\epsilon^{\prime},\alpha^{\prime})+\mathbb{N}^{n+1}$ were to hold, 
  it would imply the two inclusions 
  \begin{equation}\label{1}
    \epsilon(\beta)\in\epsilon^{\prime}+\mathbb{N}
  \end{equation}
  and 
  \begin{equation}\label{2}
    \beta\in\alpha^{\prime}+\mathbb{N}^n\,.
  \end{equation}
  Condition (ii) guarantees that not both \eqref{1} and \eqref{2} are satisfied at the same time,
  since $\epsilon(\beta)=|g|-|\beta|$ and $\epsilon^{\prime}=|g^{\prime}|-|\alpha^{\prime}|$.
  Necessity of (ii) is proved.

  Sufficiency of (i) and (ii) also follow, as we have seen that conditions (i) and (ii) translate the fact that 
  $\mathscr{G}_{\rm deglex}(\mathscr{A}_{1})$ is a reduced Gr\"obner basis into conditions on the coefficients $c_{\beta}$
  occurring in \eqref{g}. 
  
  Finally, uniqueness of $g$ satisfying (i) and (ii) follows from uniqueness of the reduced Gr\"obner basis. 
  More precisely, $\mathscr{G}_{\rm deglex}(\mathscr{A}_{1})$ contains a unique element $h$ such that 
  $\mathscr{E}_{\rm deglex}(h)=(\epsilon,\alpha)$; 
  this polynomial arises as homogenisation of $g$ satisfying (i) and (ii), which is therefore also unique. 
\end{proof}

Note that existence of a polynomial satisfying (i) such that also (ii) holds is a priori 
(i.e., by elementary properties of the polynomial ring) not clear, nor is its uniqueness.
This follows only from existence and uniqueness of the reduced Gr\"obner basis.

The necessary and sufficient criteria of Proposition \ref{onetwo} for $\alpha\in C_{\rm lex}(A_{1})$ 
to lie in $\pi(\mathscr{C}_{\rm deglex}(\mathscr{A}_{1}))$
involve only those $\alpha^{\prime}\in C_{\rm lex}(A_{1})$ for which $\alpha^{\prime}<_{\rm lex}\alpha$. 
This enables us to find an algorithm which iteratively computes the elements of 
$\mathscr{G}_{\rm deglex}(\mathscr{A}_{1})$. 
Let us first give a detailed description of the algorithm and only afterwards present its pseudocode. 

In the algorithm, we will pursue induction over the elements of a finite subset 
$C^{\prime}$ of $C_{\rm lex}(A_{1})$, 
for which we wish to decide if they belong to $\pi(\mathscr{C}_{\rm deglex}(\mathscr{A}_{1}))$, and if yes, 
construct the corresponding element of $\mathscr{G}_{\rm deglex}(\mathscr{A}_{1})$. 
The first step is to find an appropriate set $C^{\prime}$;
this should be finite and contain all $\pi(\mathscr{E}_{\rm deglex}(h))$, 
$h\in\mathscr{G}_{\rm deglex}(\mathscr{A}_{1})$. 

\begin{lmm}\label{m}
  Assume $m$ is a natural number strictly larger than $|\beta|+1$, for all $\beta$ in $D_{\rm lex}(A)$.
  Then for all $h\in\mathscr{G}_{\rm deglex}(\mathscr{A}_{1})$, we have 
  \begin{equation*}
    \mathscr{E}_{\rm deglex}(h)\in\mathscr{C}^{\prime}=
    \{\gamma\in\mathbb{N}^{n+1};\,|\gamma|\leq m\}\,,
  \end{equation*} 
  and therefore, 
  \begin{equation*}
    \pi(\mathscr{B}_{\rm deglex}(\mathscr{A}_{1}))\subset C^{\prime}
    =\{\alpha\in\mathbb{N}^{n};\,|\alpha|\leq m\}\,.
  \end{equation*}
\end{lmm}

\begin{proof}
  We start with the second assertion. 
  Let $\alpha$ be an element of $C_{\rm lex}(A_{1})$ which does not lie in the above defined $C^{\prime}$. 
  Consider $\alpha^{\prime}=\alpha-e_{j}$, where $j$ is an arbitrary element of $\{2,\ldots,n+1\}$.
  Then $\alpha^{\prime}$ also lies in $C_{\rm lex}(A_{1})$. 
  Indeed, otherwise we would have $\alpha-e_{j}\in D_{\rm lex}(A_{1})$, 
  hence in particular $|\alpha-e_{j}|<m-1$ by definition of $m$, and, on the other hand,
  $|\alpha-e_{j}|=|\alpha|-1>m-1$, as $\alpha\notin C^{\prime}$. 
  Now by definition of $C^{\prime}$, for all $\beta\in D_{\rm lex}(A_{1})$, we have $|\alpha|>|\beta|+1$, 
  hence $|\alpha^{\prime}|>|\beta|$. 
  Therefore, the total degree of polynomial $f_{\alpha}$ equals $|\alpha|$, 
  and the total degree of polynomial $f_{\alpha^{\prime}}$ equals $|\alpha^{\prime}|$. 
  It follows that $\mathscr{E}_{\rm deglex}(\mathscr{H}(f_{\alpha}))=(0,\alpha)$, 
  and $\mathscr{E}_{\rm deglex}(\mathscr{H}(f_{\alpha^{\prime}}))=(0,\alpha^{\prime})$. 
  In particular, $(0,\alpha)$ lies in is $(0,\alpha^\prime)+\mathbb{N}^{n+1}$, 
  which implies that $(0,\alpha)$ not an element of $\mathscr{B}_{\rm deglex}(\mathscr{A}_{1})$. 
  Therefore, also no element of the form $(\epsilon,\alpha)$, where $\epsilon\in\mathbb{N}$, 
  lies in $\mathscr{B}_{\rm deglex}(\mathscr{A}_{1})$, and the second assertion is shown.
  
  Now we prove the first assertion. We will have to make a clear distinction between the leading term of certain polynomials
  w.r.t. $\leq_{\rm lex}$ and w.r.t. $\preceq_{\rm deglex}$. This will always be pointed out explicitly.
  Take an element $h$ of $\mathscr{G}_{\rm deglex}(\mathscr{A})$, write $\gamma=\mathscr{E}_{\rm deglex}(h)$,
  $\alpha=\pi(\gamma)$, and $\epsilon$ for the first component of $\gamma$; 
  then we have to show that $\epsilon$ does not exceed $m-|\alpha|$.
  Consider the polynomial $f_{\alpha}\in I(A_{1})$. By what we have already shown, $\alpha$ lies in $C^{\prime}$, 
  and all lexicographcally nonleading exponents $\beta$ of $f_{\alpha}$ lie in $D_{\rm lex}(A_{1})$ by definition. 
  Therefore, the graded lexicographic leading exponent of $f_{\alpha}$ is 
  \begin{itemize}
    \item either $\alpha$, 
    \item or $\beta$, for some $\beta\in D_{\rm lex}(A_{1})$. 
  \end{itemize}

  In the first case, the graded lexicographically leading exponent of $\mathscr{H}(f_{\alpha})$ is $(0,\alpha)$. 
  As $\mathscr{H}(f_{\alpha})$ lies in $\mathscr{I}(\mathscr{A}_{1})$, 
  the exponent $\gamma$ we started with must also equal $(0,\alpha)$. 
  By the second assertion of the Lemma, $\alpha$ lies in $C^{\prime}$, 
  hence $\gamma$ lies in $\mathscr{C}^\prime$. 

  In the second case, the graded lexicographically leading term of $\mathscr{H}(f_{\alpha})$ 
  is $(\epsilon^\prime,\alpha)$, where $\epsilon^\prime=|\beta|-|\alpha|$. 
  However, definition of $m$ implies, in particular, that $D_{\rm lex}(A_{1})$ is contained in 
  $C^{\prime}$, hence $|\beta|\leq m$ for all $\beta$, hence $|(\epsilon^\prime,\alpha)|\leq m$. 
  It follows that $(\epsilon^\prime,\alpha)$, the graded lexicographically leading exponent of 
  $\mathscr{H}(f_{\alpha})$, lies in $\mathscr{C}^{\prime}$.
  This implies that for the exponent $\gamma$ we started with, 
  we must have $\epsilon\leq\epsilon^\prime$. It follows that $|\gamma|\leq|(\epsilon^\prime,\alpha)|\leq m$, 
  hence $\gamma\in\mathscr{C}^{\prime}$, as claimed.
\end{proof}

Now we fix an element $\alpha\in C^{\prime}$. 
We assume for all $\alpha^{\prime}\in\pi(\mathscr{B}_{\rm deglex}(\mathscr{A}_{1}))$ such that 
$\alpha^{\prime}<_{\rm lex}\alpha$, 
we are given the corresponding element of $\mathscr{G}_{\rm deglex}(\mathscr{A}_{1})$, 
i.e., a polynomial $h^{\prime}\in\mathscr{G}_{\rm deglex}(\mathscr{A}_{1})$ such that
$\pi(\mathscr{E}_{\rm deglex}(h^\prime))=\alpha^\prime$. 
We wish to decide whether or not $\alpha$ leads to a polynomial $g$ of the form \eqref{g} 
such that $h=\mathscr{H}(g)\in\mathscr{G}_{\rm deglex}(\mathscr{A}_{1})$. 

Let us assume for the moment that our given $\alpha$ does lead to an element 
$h=\mathscr{H}(g)$ of $\mathscr{G}_{\rm deglex}(\mathscr{A}_{1})$. 
A priori we do now know the total degree of $g$; we can only say that $|g|\geq|\alpha|$. 
Therefore, we will have
\begin{equation*}
  |g|=|\alpha|+r\,,
\end{equation*}
for some $r\in\mathbb{N}$ yet to be determined. Condition (ii) of Proposition \ref{onetwo} tells us that 
in the sum \eqref{g}, not all $\beta\in C_{\rm lex}(A)$ such that $\beta<_{\rm lex}\alpha$ will lead to a nonzero $c_{\beta}$;
we can exclude some of them. More precisely,
\begin{equation}\label{gy}
  g=f_{\alpha}+\sum_{\beta\in Y_r}c_{\beta}X^{\beta}\,,
\end{equation}
where
\begin{equation*}
  \begin{split}
    Y_r=\{&\beta\in C_{\rm lex}(A);\,\beta<_{\rm lex}\alpha,\\
    &\text{ and for all }\alpha^{\prime}\in\pi(\mathscr{G}_{\rm deglex}(\mathscr{A}))
    \text{ such that }\beta\in\alpha^{\prime}+\mathbb{N}^{n},\\
    &\text{ we have }|\alpha|+r+|\alpha^{\prime}|-|g^{\prime}|<|\beta|\leq|\alpha|+r\}\,.
  \end{split}
\end{equation*}
Indeed, inequality $|\alpha|+r+|\alpha^{\prime}|-|g^{\prime}|<|\beta|$ is imposed by Proposition \ref{onetwo}, (ii), 
and inequality $|\beta|\leq|\alpha|+r$ is imposed by the fact that the total degree of none of the $f_{\beta}$ 
involved in $g$ may exceed $|\alpha|+r$. 

Let us now find a linear system of equations for the coefficients $c_{\beta}$ in equation \eqref{gy}, 
which uniquely determines the $c_{\beta}$. 
For this, we need all coefficients of all $f_{\beta}$. For $\beta\in C_{\rm lex}(A_{1})$, let us write 
\begin{equation*}
  f_{\beta}=X^{\beta}+\sum_{\substack{\beta^{\prime}\in D_{\rm lex}(A_{1})\\\beta^{\prime}<_{\rm lex}\beta}}
  d_{\beta,\beta^{\prime}}X^{\beta^{\prime}}\,,
\end{equation*}
or, in a more compact way, 
\begin{equation*}
  f_{\beta}=\sum_{\beta^{\prime}\in\mathbb{N}^{n}}d_{\beta,\beta^{\prime}}X^{\beta^{\prime}}\,,
\end{equation*}
where $d_{\beta,\beta}=1$, 
and $d_{\beta,\beta^{\prime}}=0$ whenever $\beta^{\prime}\notin D_{\rm lex}(A_{1})\cup\{\beta\}$ 
or $\beta^{\prime}>\beta$. 
Therefore, equation \eqref{gy} says that
\begin{equation}\label{gd}
  g=\sum_{\beta^{\prime}\in\mathbb{N}^{n}}(d_{\alpha,\beta^{\prime}}
  +\sum_{\beta\in Y_{r}}d_{\beta,\beta^{\prime}}c_{\beta})X^{\beta^{\prime}}\,.
\end{equation}
Our assumption is that the total degree of $g$ equals $|\alpha|+r$. 
This means that we have to make sure that the coefficients of all monomials $X^{\beta^{\prime}}$ occurring in 
the above equation for which $|\beta^{\prime}|>|\alpha|+r$, are in fact zero. 
This is expressed by the following system of linear equations in the variables $c_{\beta}$, for $\beta\in Y_{r}$. 
\begin{equation}\label{sys}
  d_{\alpha,\beta^{\prime}}+\sum_{\beta\in Y_{r}}d_{\beta,\beta^{\prime}}c_{\beta}=0
  \text{ for all }\beta^{\prime}\text{ such that }|\beta^{\prime}|>|\alpha|+r\,.
\end{equation}

If $\alpha$ indeed leads to an element $h=\mathscr{H}(g)$ of $\mathscr{G}_{\rm deglex}(\mathscr{A}_{1})$, 
and if we picked the right $r$ such that $|h|=|g|=|\alpha|+r$, 
then system \eqref{sys} has a unique by Proposition \ref{onetwo}.
However, we still do not know 
\begin{itemize}
  \item the total degree, $|\alpha|+r$, of $g$, 
  \item nor even if $\alpha$ leads to $h=\mathscr{H}(g)\in\mathscr{G}_{\rm deglex}(\mathscr{A}_{1})$ or not. 
\end{itemize}

But both these questions can be decided by the following iterative process. 
We start by assuming that $r=0$. This leads to a set $Y_{0}$, as defined above, 
and accordingly, to a system of equations \eqref{sys}, 
where $r=0$. If this system has a solution, it will by Proposition \ref{onetwo} be unique, and by that same Proposition, 
$h=\mathscr{H}(g)$ will lie in $\mathscr{G}_{\rm deglex}(\mathscr{A}_{1})$. 
If system \eqref{sys} does not have a solution, we try $r=1$. Analogously as before, 
we compute set $Y_{1}$ and decide if system of equations \eqref{sys}, for $r=1$, has a (unique) solution. 
If yes, $h=\mathscr{H}(g)$ will lie in $\mathscr{G}_{\rm deglex}(\mathscr{A}_{1})$. If no, try $r=2$, etc. 
If by this process, we find an $r$ such that system of equations \eqref{sys} has a solution, 
we take the smallest such $r$ and deduce that $h=\mathscr{H}(g)$ lies in 
$\mathscr{G}_{\rm deglex}(\mathscr{A}_{1})$, and that $|g|=|\alpha|+r$. 

It may happen that for all $r\geq 0$, system \eqref{sys} does not have a solution. 
By Proposition \ref{onetwo}, this means that $\alpha$ does not lie in $\pi(\mathscr{B}_{\rm deglex}(\mathscr{A}_{1}))$, 
i.e. there exists no $g$ as in \eqref{g} whose homogenisation is an element of the 
Gr\"obner basis $\mathscr{G}_{\rm deglex}(\mathscr{A}_{1})$.
Fortunately, we can decide if this is the case by considering only a finite number of $r$, namely, those $r$ for which
$0\leq r\leq|\alpha|-m$. This follows from Lemma \ref{m}.

Before summarising our iterative construcion of $\mathscr{G}_{\rm deglex}(\mathscr{A}_{1})$ in Algorithm \ref{algone}, 
we make one last remark concerning the set $C^{\prime}$, which is our set of candidates for elements of 
$\pi(\mathscr{B}_{\rm deglex}(\mathscr{A}_{1}))$. After having examined $\alpha$ in the process described above, 
we can of course delete $\alpha$ from the set of candidates, $C^{\prime}$ -- no matter if $\alpha$ is 
contained in $\pi(\mathscr{B}_{\rm deglex}(\mathscr{A}_{1}))$ or not.
However, if $(0,\alpha)$ lies in $\mathscr{B}_{\rm deglex}(\mathscr{A}_{1})$
(i.e. $r=0$ in the notation used above), no element of $\alpha+\mathbb{N}^{n}$ will lie in 
$\pi(\mathscr{B}_{\rm deglex}(\mathscr{A}_{1}))$. Hence in this case, we can delete $\alpha+\mathbb{N}^{n}$ from
the set of candidates, $C^{\prime}$. 
Doing so, we shrink the finite set $C^{\prime}$ in each step of the process. 
This suggests to take $C^{\prime}$ as control variable, as is pursued in Algorithm \ref{algone}.

\begin{algorithm}
  \caption{Calculate $\mathscr{G}_{\rm deglex}(\mathscr{A}_{1})$}\label{algone}
  \begin{algorithmic}
    \STATE $C^{\prime}:=\{0,\ldots,m\}\times\ldots\times\{0,\ldots,m\}\subset\mathbb{N}^n$
    \STATE $\mathscr{G}_{\rm deglex}(\mathscr{A}_{1}):=\emptyset$
    \WHILE{$C^{\prime}\neq\emptyset$}
      \STATE $\alpha:=$ the lexicographic minimum of $C^{\prime}$
      \STATE $r:=0$
      \WHILE{$r\leq m-|\alpha|$}
        \STATE check if \eqref{sys} has a solution
        \IF{yes}
          \STATE $g:=$ the polynomial defined in \eqref{gd}
          \STATE $\mathscr{G}_{\rm deglex}(\mathscr{A}_{1})
            :=\mathscr{G}_{\rm deglex}(\mathscr{A}_{1})\cup\{\mathscr{H}(g)\}$
          \STATE $r:=m-|\alpha|+2$
        \ELSE
          \STATE $r:=r+1$
        \ENDIF
      \ENDWHILE
      \IF{$r=m-|\alpha|+1$}
        \STATE $C^{\prime}:=C^{\prime}-\{\alpha\}$
      \ELSE
        \STATE $C^{\prime}:=C^{\prime}-(\alpha+\mathbb{N}^n)$
      \ENDIF
    \ENDWHILE
  \end{algorithmic}
\end{algorithm}

\section{Axes in $\mathscr{D}(\mathscr{A})$}\label{allaxes}

Let us now consider the following sets.
\begin{equation*}
  \begin{split}
    A_{1}&=\mathscr{A}\cap\{X_{1}=1\}\,,\\
    A_{2}&=\mathscr{A}\cap\{X_{1}=0\}\cap\{X_{2}=1\}\,,\\
    \vdots&\\
    A_{j}&=\mathscr{A}\cap\{X_{1}=0\}\cap\ldots\cap\{X_{j-1}=0\}\cap\{X_{j}=1\}\,,\\
    \vdots&\\
    A_{n+1}&=\mathscr{A}\cap\{X_{1}=0\}\cap\ldots\cap\{X_{n}=0\}\cap\{X_{n+1}=1\}\,;\\
  \end{split}
\end{equation*}
We denote by $\langle A_{j}\rangle$ the union of all lines passing through elements of $A_{j}$;
then clearly 
\begin{equation*}
  \mathscr{A}=\langle A_{1}\rangle\coprod\ldots\coprod\langle A_{n+1}\rangle\,,
\end{equation*}
and, with the notation introduced in the Introduction, 
\begin{equation*}
  \begin{split}
    \mathscr{A}_{0}&=\langle A_{2}\rangle\coprod\ldots\coprod\langle A_{n+1}\rangle\,,\\
    \mathscr{A}_{1}&=\langle A_{1}\rangle\,.
  \end{split}
\end{equation*}

In Section \ref{axes}, we studied the set $\mathscr{A}_{1}$, where we proved Theorem \ref{thma1}, 
which was a statement about the axes contained in $\mathscr{D}(\mathscr{A}_{1})$.
In this Section, we will study the entire set $\mathscr{A}$, 
and prove a generalisation of the statement of Theorem \ref{thma1}, 
concering the axes contained in $\mathscr{D}(\mathscr{A})$.
The technique will be induction over $n$, the dimension of the projective space we are working in. 
Note that $\mathscr{A}_{0}$ is contained in $\mathbb{P}^{n-1}$;
this observation is going to be of fundamental impontance of the induction step. 

The above defined set $A_{j}$ is naturally embedded into $\mathbb{A}^{n+1-j}$ by coordinates $X_{j+1},\ldots,X_{n+1}$. 
This yields, in particular, $D(A_{j})\subset\mathbb{N}^{n+1-j}$. Consider the embedding
\begin{equation*}
  \begin{split}
    i_{n+1-j}:\mathbb{N}^{n+1-j}&\hookrightarrow\mathbb{N}^{n+1}\\
    \sigma=(\sigma_{j+1},\ldots,\sigma_{n+1})&\mapsto(0,\sigma)=(0,\ldots,0,\sigma_{j+1},\ldots,\sigma_{n+1})\,,
  \end{split}
\end{equation*}
where we fill up the $(n+1)$-tuple by $j$ times the number $0$. 
The embedding $A_{j}\hookrightarrow\mathbb{A}^{n+1-j}$ yields also an embedding 
$\langle A_{j}\rangle\hookrightarrow\mathbb{P}^{n-j}$.
Note that all elements $\ell$ of $\langle A_{j}\rangle$ lie in the finite part of $\mathbb{P}^{n-j}$, i.e., 
the intersection of $\ell$ with the set $\{X_{j+1}=1\}$ is nontrivial. 
Hence by Theorem \ref{thma1}, all affine axes contained in $\mathscr{D}(\langle A_{j}\rangle)$, 
when embedded into $\mathbb{N}^{n+1}$ by $i_{n+1-j}$, are parallel to $\mathbb{N}e_{j}$. 
Furthermore, the number of affine axes in $\mathscr{D}(\langle A_{j}\rangle)$ equals $\#D(A_{j})$.

\begin{thm}\label{thma}
  For all $j\in\{1,\ldots,n+1\}$, the set $\mathscr{D}(\mathscr{A})$ contains precisely 
  $\#D(A_{j})$ axes parallel to $\mathbb{N}e_{j}$. 
  These are the images by $i_{n+1-j}$ of the axes contained in $\mathscr{D}(\langle A_{j}\rangle)$.
\end{thm}

\begin{proof}
  Let $h$ be an element of $\mathscr{I}(\mathscr{A})$. Then from $h\in\mathscr{I}(\mathscr{A}_{0})$, 
  we conclude that $\mathscr{E}(h)\in\mathscr{C}(\mathscr{A}_{0})$, and from $h\in\mathscr{I}(\mathscr{A}_{1})$, 
  we conclude that $\mathscr{E}(h)\in\mathscr{C}(\mathscr{A}_{1})$. This shows that 
  $\mathscr{C}(\mathscr{A})\subset\mathscr{C}(\mathscr{A}_{0})\cap\mathscr{C}(\mathscr{A}_{1})$, 
  or, in terms of complements,
  \begin{equation}\label{supset}
    \mathscr{D}(\mathscr{A})\supset\mathscr{D}(\mathscr{A}_{0})\cup\mathscr{D}(\mathscr{A}_{1})\,.
  \end{equation}
  By induction over $n$, we may assume that for all $j\in\{2,\ldots,n+1\}$, 
  the set $\mathscr{D}(\mathscr{A}_{0})$ contains precisely $\#D(A_{j})$ axes parallel to $\mathbb{N}e_{j}$, 
  namely, the images by $i_{n+1-j}$ of the axes contained in $\mathscr{D}(\langle A_{j}\rangle)$.
  Furthermore, by Theorem \ref{thma1}, 
  there are precisely $\#D(A_{1})$ axes contained in $\mathscr{D}(\mathscr{A}_{1})$, all parallel to $\mathbb{N}e_{1}$, 
  namely, the axes contained in $\mathscr{D}(\langle A_{1}\rangle)$. (Note that $i_{n+1}$ is the identical map.)
  Therefore, by \eqref{supset}, we have to prove that for all $j\in\{1,\ldots,n+1\}$, 
  the set $\mathscr{D}(\mathscr{A})$ contains not more 
  axes parallel to $\mathbb{N}e_{j}$ than those we already found. 
  
  First, we show the assertion for $j=1$. 
  Consider the ideal $\mathscr{I}(\mathscr{A}_{1})$, and for all $\gamma\in\mathscr{B}(\mathscr{A}_{1})$, 
  the element $f_{\gamma}$ of the Gr\"obner basis $\mathscr{G}(\mathscr{A}_{1})$. 
  The product $X_{1}f_{\gamma}$ lies in $\mathscr{I}(\mathscr{A})$, and $\mathscr{E}(X_{1}f_{\gamma})=(1,\gamma)$.
  This shows that $\mathscr{D}(\mathscr{A})$ contains not more axes parallel to $\mathbb{N}e_{1}$ than those
  contained in $\mathscr{D}(\langle A_{1}\rangle)$.
  
  Now we fix a $j\in\{2,\ldots,n+1\}$. By Theorem \ref{thma1}, 
  the only axes contained in $\mathscr{D}(\mathscr{A}_{1})$ are parallel to $\mathbb{N}e_{1}$, 
  hence for a $\mu\in\mathbb{N}$ which is sufficiently large, we have $\mu e_{j}\in\mathscr{C}(\mathscr{A}_{1})$. 
  Lemma \ref{fund}, applied to the ideal $\mathscr{I}(\mathscr{A}_{1})$, 
  provides a polynomial $f^{(1)}_{\mu e_{j}}\in\mathscr{I}(\mathscr{A}_{1})$
  such that $\mathscr{E}(f^{(1)}_{\mu e_{j}})=\mu e_{j}$.
  Furthermore, we fix an arbitrary $\gamma\in\mathscr{C}(\mathscr{A}_{0})$. 
  Lemma \ref{fund}, applied to the ideal $\mathscr{I}(\mathscr{A}_{0})$, 
  provides a polynomial $f^{(0)}_{\gamma}\in\mathscr{I}(\mathscr{A}_{0})$ 
  such that $\mathscr{E}(f^{(0)}_{\gamma})=\gamma$. 
  It follows that $f^{(1)}_{\mu e_{j}}f^{(0)}_{\gamma}\in\mathscr{I}(\mathscr{A})$, 
  and $\mathscr{E}(f^{(1)}_{\mu e_{j}}f^{(0)}_{\gamma})
  =\mathscr{E}(f^{(1)}_{\mu e_{j}})+\mathscr{E}(f^{(1)}_{\gamma})=\mu e_{j}+\gamma$. 
  Therefore, for all $\gamma\in\mathscr{C}(\mathscr{A}_{0})$, 
  the axis through $\gamma$ parallel to $\mathbb{N}e_{j}$ is not contained in $\mathscr{D}(\mathscr{A})$. 
  Conversely, each axis parallel to $\mathbb{N}e_{j}$ which is contained in $\mathscr{D}(\mathscr{A})$
  must be contained in $\mathscr{D}(\mathscr{A}_{0})$. 
  This shows that $\mathscr{D}(\mathscr{A})$ contains not more axes parallel to $\mathbb{N}e_{j}$ than those
  we already found at the beginning.
\end{proof}

\begin{cor}
  The axes contained in $\mathscr{D}_{\rm deglex}(\mathscr{A})$ are given by
  \begin{equation*}
    i_{n+1-j}(\sigma)+\mathbb{N}e_{j}\,,
  \end{equation*}
  for $\sigma\in D_{\rm lex}(A_{j})$ and $j\in\{1,\ldots,n+1\}$. 
\end{cor}

\section{The ideal of $\mathscr{A}$}\label{merge}

Similarly as in Section \ref{stepone}, we now present an algorithm for the construction of the Gr\"obner basis
$\mathscr{G}(\mathscr{A})$ of $\mathscr{I}(\mathscr{A})$. 
The starting point is the Gr\"obner bases of $\langle A_{j}\rangle$, for $j=1,\ldots,n+1$. 
(It is justified to assume these Gr\"obner bases given, as we we studied these bases in Section \ref{stepone}
for the term order $\leq_{\rm lex}$, and the Gr\"obner bases for all other term orders can be computed therefrom.) 
Our technique is the same as in Section \ref{allaxes} -- 
we build up $\mathscr{G}(\mathscr{A})$ by induction over the dimension. 
\begin{itemize}
  \item We start with $\langle A_{n+1}\rangle$, whose Gr\"obner basis is known. 
  \item Then we merge the Gr\"obner bases of $\langle A_{n}\rangle$, 
    a subset of the finite part of $\mathbb{P}^{1}$, 
    with the Gr\"obner basis of $\langle A_{n+1}\rangle$, 
    which we consider a subset of the infinite part of $\mathbb{P}^{1}$. 
  \item Then we merge the Gr\"obner bases of $\langle A_{n-1}\rangle$, 
    a subset of the finite part of $\mathbb{P}^{2}$, 
    with the Gr\"obner basis of $\langle A_{n}\rangle\coprod\langle A_{n+1}\rangle$, 
    which we consider a subset of the infinite part of $\mathbb{P}^{2}$. 
  \item In this way we proceed, util finally, 
    we merge the Gr\"obner bases of $\langle A_{1}\rangle$, a subset of the finite part of $\mathbb{P}^{n}$, 
    with the Gr\"obner basis of $\langle A_{2}\rangle\coprod\ldots\coprod\langle A_{n+1}\rangle$, 
    considered as a subset of the infinite part of $\mathbb{P}^{n}$. 
\end{itemize}
All steps are essentially one and the same thing -- it all comes down to the very last step, i.e., 
merging the Gr\"obner bases of $\mathscr{A}_{1}$, a subset of the finite part of $\mathbb{P}^{n}$, 
with the Gr\"obner basis of $\mathscr{A}_{0}$, considered as a subset of the infinite part of $\mathbb{P}^{n}$. 
Hence, the input of our algorithm for the construction of $\mathscr{G}(\mathscr{A})$ is
\begin{itemize}
  \item the Gr\"obner basis $\mathscr{G}(\mathscr{A}_{0})$ of $\mathscr{I}(\mathscr{A}_{0})$, and
  \item the Gr\"obner basis $\mathscr{G}(\mathscr{A}_{1})$ of $\mathscr{I}(\mathscr{A}_{1})$.
\end{itemize}
Given these data, Lemma \ref{fund} provides
\begin{itemize}
  \item for all $\gamma\in\mathscr{C}(\mathscr{A}_{0})$, a polynomial $f_{\gamma}^{(0)}\in\mathscr{I}(\mathscr{A}_{0})$
  such that $\mathscr{E}(f_{\gamma}^{(0)})=\gamma$, and all nonleading terms of $f_{\gamma}^{(0)}$ lie in 
  $\mathscr{D}{(\mathscr{A}_{0})}$, and
  \item for all $\gamma\in\mathscr{C}(\mathscr{A}_{1})$, a polynomial $f_{\gamma}^{(1)}\in\mathscr{I}(\mathscr{A}_{1})$
  such that $\mathscr{E}(f_{\gamma}^{(1)})=\gamma$, and all nonleading terms of $f_{\gamma}^{(1)}$ lie in 
  $\mathscr{D}{(\mathscr{A}_{1})}$. 
\end{itemize}

In the algorithm, the control variable will be a set $\mathscr{C}^{\prime}\subset\mathbb{N}^{n+1}$, 
which is the set of candidates for exponents of elements of $\mathscr{G}(\mathscr{A})$. 
As $\mathscr{C}(\mathscr{A})\subset\mathscr{C}(\mathscr{A}_{0})\cap\mathscr{C}(\mathscr{A}_{1})$, 
we may start with $\mathscr{C}^{\prime}=\mathscr{C}(\mathscr{A}_{0})\cap\mathscr{C}(\mathscr{A}_{1})$. 
Set $\mathscr{C}^{\prime}$ will be shrunk in each step of the algorithm.
In each step, we let $\gamma$ be the minimal (w.r.t. $\preceq$) element of $\mathscr{C}^{\prime}$.
We will decide if $\gamma\in\mathscr{B}(\mathscr{A})$ or not. 
If yes, we will compute the corresonding $f_{\gamma}\in\mathscr{G}(\mathscr{A})$, 
and we will replace $\mathscr{C}^{\prime}$ by 
$\mathscr{C}^{\prime}-(\gamma+\mathbb{N}^{n+1})$. 
If no, we will replace $\mathscr{C}^{\prime}$ by $\mathscr{C}^{\prime}-\{\gamma\}$. 
In the next step of the algorithm, we pass to the minimal (w.r.t. $\preceq$) element of the new 
(i.e., smaller) set $\mathscr{C}^{\prime}$. 

The question whether or not $\gamma$ lies in $\mathscr{B}(\mathscr{A})$ is decided as follows. 
Firstly, of course, $\gamma$ lies in $\mathscr{C}(\mathscr{A})$ if, and only if, 
there exists $f_{\gamma}\in\mathscr{I}(\mathscr{A})$ such that $\mathscr{E}(f_{\gamma})=\gamma$. 
Polynomial $f_{\gamma}$, if it exists, in particular lies in $\mathscr{I}(\mathscr{A}_{0})$, 
and is therefore of the form 
\begin{equation}\label{f1}
  f_{\gamma}=f_{\gamma}^{(0)}+
  \sum_{\substack{\delta\in\mathscr{C}(\mathscr{A}_{0})\\ |\delta|=|\gamma|\\ \delta\prec\gamma}}
  c_{\delta}f_{\delta}^{(0)}\,,
\end{equation}
for some $c_{\delta}\in k$. Furthermore, polynomial $f_{\gamma}$, if it exists, 
also lies in $\mathscr{I}(\mathscr{A}_{1})$, and is therefore also of the form 
\begin{equation}\label{f2}
  f_{\gamma}=f_{\gamma}^{(1)}+
  \sum_{\substack{\eta\in\mathscr{C}(\mathscr{A}_{1})\\ |\eta|=|\gamma|\\ 
  \eta\prec\gamma}}d_{\eta}f_{\eta}^{(1)}\,,
\end{equation}
for some $d_{\eta}\in k$. 
(Note that $\gamma$ lies in both $\mathscr{C}(\mathscr{A}_{0})$ and $\mathscr{C}(\mathscr{A}_{1})$, 
as we have seen above.) Therefore, $\gamma$ lies in $\mathscr{C}(\mathscr{A})$ if, and only if, equation
\begin{equation}\label{twosums}
  f_{\gamma}^{(0)}-f_{\gamma}^{(1)}=
  \sum_{\substack{\eta\in\mathscr{C}(\mathscr{A}_{1})\\ |\eta|=|\gamma|\\ 
  \eta\prec\gamma}}d_{\eta}f_{\eta}^{(1)}
  -\sum_{\substack{\delta\in\mathscr{C}(\mathscr{A}_{0})\\ |\delta|=|\gamma|\\ \delta\prec\gamma}}
  c_{\delta}f_{\delta}^{(0)}
\end{equation}
(a linear equation in variables $c_{\delta}$ and $d_{\eta}$), has a solution. 

However, if \eqref{twosums} has a solution, the solution will in general not be unique. 
This is due to the fact that there may be many different polynomials in $\mathscr{I}(\mathscr{A})$
with the same leading term. For getting uniqueness of $f_{\gamma}$, 
we have to impose a restriction also on the nonleading terms of $f_{\gamma}$.
By Lemma \ref{fund}, $f_{\gamma}$ will be unique if we also require that no nonleading term of $f_{\gamma}$ 
is divided by any $X^\xi$, where $\xi$ is an element of $\mathscr{B}(\mathscr{A})$ such that $\xi\prec\gamma$. 
(Note by the recursive character of our algorithm, 
when we are at the point where we check $\gamma\in\mathscr{C}^\prime$, 
we know all $\xi\in\mathscr{B}(\mathscr{A})$ such that $\xi\prec\gamma$, along with the corresponding $f_{\xi}$.)
In terms of a system of linear equations, this condition reads as follows.
\begin{equation}\label{additionally}
  \begin{split}
    \forall\xi\in\mathscr{B}(\mathscr{A})\text{ s.t. }\xi\prec\gamma\,,
    \text{ if }&\delta\in\xi+\mathbb{N}^{n+1}\,,\text{ then }c_{\delta}=0\,,\\
    \text{ and if }&\eta\in\xi+\mathbb{N}^{n+1}\,,\text{ then }d_{\eta}=0\,.
  \end{split}
\end{equation}
Indeed, in equations \eqref{f1} and \eqref{f2} for $f_{\gamma}$, only the leading terms of $f_{\delta}$ (resp. $f_{\eta}$)
might be divided by some $X^\xi$ as above, as all nonleading exponents of $f_{\delta}$ (resp. $f_{\eta}$)
lie in $\mathscr{D}(\mathscr{A})$. 

Thus $\gamma$ lies in $\mathscr{C}(\mathscr{A})$ if, and only if, 
the system of equations defined by \eqref{twosums} and \eqref{additionally} has a solution. 
Now we can prove an even stronger assertion holds, 
taking advantage of the recursive definition of $\mathscr{C}^\prime$. 

\begin{lmm}
  Let $\gamma$ be the minimum w.r.t. $\preceq$ of the recursively defined set $\mathscr{C}^{\prime}$.
  Then $\gamma$ lies in $\mathscr{B}(\mathscr{A})$ if, and only if, 
  the system of equations (in variables $c_{\delta}$ and $d_{\eta}$) 
  defined by \eqref{twosums} and \eqref{additionally} has a solution. 
  In this case, $f_{\gamma}$ lies in $\mathscr{G}(\mathscr{A})$.
\end{lmm}

\begin{proof}
  It is clear that if $f_{\gamma}$ is an element of $\mathscr{G}(\mathscr{A})$
  such that $\mathscr{M}(f_{\gamma})=X^\gamma$, 
  the system defined by \eqref{twosums} and \eqref{additionally} has a solution. 
  This is independent of $\gamma$ being the minimum of $\mathscr{C}^{\prime}$.
  
  Let us show that if $\gamma$ is as in the Lemma, 
  and if systems \eqref{twosums} and \eqref{additionally} have a solution, 
  then $\gamma\in\mathscr{B}(\mathscr{A})$, and $f_{\gamma}\in\mathscr{G}(\mathscr{A})$.
  We verify this by induction over the recursively defined set $\mathscr{C}^{\prime}$.

  At the basis of the induction, set $\mathscr{C}^{\prime}$ has the property that no element of 
  $\mathscr{C}^{\prime}$ lies in $\mathscr{B}(\mathscr{A})$. 
  This remains valid as long as we find only such $\gamma\in\mathscr{C}^{\prime}$ for which
  \eqref{twosums} and \eqref{additionally} have no solution. 
  At some point, we will reach the first $\gamma\in\mathscr{C}^{\prime}$ 
  for which there exists a solution to \eqref{twosums} and \eqref{additionally}.
  Then we have $\gamma\in\mathscr{B}(\mathscr{A})$. 
  Indeed, we know that $\gamma\in\mathscr{C}(\mathscr{A})$, 
  hence there exists $\gamma^{\prime}\in\mathscr{B}(\mathscr{A})$ is such that 
  $\gamma\in\gamma^{\prime}+\mathbb{N}^{n+1}$. Then, on the one hand, 
  $\gamma^{\prime}$ leads to a solution of a system analogous to \eqref{twosums} and \eqref{additionally}, 
  where each $\gamma$ is replaced by $\gamma^{\prime}$, and on the other hand, 
  we have $\gamma^{\prime}\preceq\gamma$. 
  By minimality of $\gamma$, we conclude that $\gamma^{\prime}=\gamma$, 
  hence $\gamma\in\mathscr{B}(\mathscr{A})$.
  
  By construction, all nonleading terms of $f_{\gamma}$ lie in $\mathscr{D}(\mathscr{A})$, 
  hence $f_{\gamma}$ is an element of the reduced Gr\"obner basis, $\mathscr{G}(\mathscr{A})$.
  
  At any later stage of the induction, we consider the minimum $\gamma$ of $\mathscr{C}^{\prime}$,
  and decide if equations \eqref{twosums} and \eqref{additionally} have a solution or not. 
  At this stage, the current $\mathscr{C}^{\prime}$ contains no element of any set
  $\gamma^\prime+\mathbb{N}^{n+1}$, for all $\gamma^\prime\in\mathscr{B}(\mathscr{A})$ 
  such that $\gamma^\prime\prec\gamma$. 
  Now if \eqref{twosums} and \eqref{additionally} have a solution, then $\gamma$ lies in $\mathscr{B}(\mathscr{A})$. 
  Indeed, we know that $\gamma\in\mathscr{C}(\mathscr{A})$, 
  hence there exists $\gamma^\prime\in\mathscr{B}(\mathscr{A})$ is such that 
  $\gamma\in\gamma^\prime+\mathbb{N}^{n+1}$, hence, in particular, $\gamma^\prime\preceq\gamma$. 
  By what we remarked about $\mathscr{C}^{\prime}$, 
  therefrom follows $\gamma^\prime=\gamma$, and $\gamma\in\mathscr{B}(\mathscr{A})$.
  
  The reason for $f_{\gamma}\in\mathscr{G}(\mathscr{A})$ is the same as in the induction basis. 
\end{proof}

Let us formulate our construction in a pseudocode, Algorithm \ref{algmerge}. 

\begin{algorithm}\label{algmerge}
  \caption{Calculate $\mathscr{G}(\mathscr{A})$}
  \begin{algorithmic}
    \STATE $\mathscr{C}^{\prime}(\mathscr{A})=\mathscr{C}(\mathscr{A}_{0})\cap\mathscr{C}(\mathscr{A}_{1})$
    \WHILE{$\mathscr{C}^{\prime}(\mathscr{A})\neq\emptyset$}
      \STATE $\gamma=$ the graded lexicographic minimum of $\mathscr{C}^{\prime}(\mathscr{A})$
      \STATE check if equations \eqref{twosums} and \eqref{additionally} 
        (in the variables $c_{\delta}$, $d_{\eta}$) have a solution
      \IF{yes}
        \STATE $\mathscr{G}(\mathscr{A})=\mathscr{G}(\mathscr{A})\cup\{f_{\gamma}\}$
        \STATE $\mathscr{C}^{\prime}(\mathscr{A})=\mathscr{C}^{\prime}(\mathscr{A})-(\gamma+\mathbb{N}^{n+1})$
      \ELSE
        \STATE $\mathscr{C}^{\prime}(\mathscr{A})=\mathscr{C}^{\prime}(\mathscr{A})-\{\gamma\}$
      \ENDIF
    \ENDWHILE
  \end{algorithmic}
\end{algorithm}

\section{Acknoledgements}\label{ack}

I wish to thank Thomas Zink for his constant encouragement for my research in commutative algebra, 
in particular, Gr\"obner bases.
Therefore, I am currently working on generalisations of Theorems \ref{thmone} through \ref{thma}, 
concerning a union of lines which do not necessarily pass through the origin, 
and also a generalisation to a finite union arbitrary linear subspaces of dimension $d$ in $\mathbb{A}^n$, 
in particular, unions of points on the Grassmannian ${\rm Gr}_{d}(n)$.
%Many thanks go to my brother Thomas Lederer for his help in questions concerning English language.

\bibliography{proj.bib}

\begin{thebibliography}{ABKR00}

\bibitem[ABKR00]{akr1}
J.~Abbott, A.~Bigatti, M.~Kreuzer, and L.~Robbiano.
\newblock Computing ideals of points.
\newblock {\em J. Symbolic Comput.}, 30(4):341--356, 2000.

\bibitem[AKR05]{akr2}
J.~Abbott, M.~Kreuzer, and L.~Robbiano.
\newblock Computing zero-dimensional schemes.
\newblock {\em J. Symbolic Comput.}, 39(1):31--49, 2005.

\bibitem[AMM03]{big}
Maria~Emilia Alonso, Maria~Grazia Marinari, and Teo Mora.
\newblock The big mother of all dualities: {M}\"oller algorithm.
\newblock {\em Comm. Algebra}, 31(2):783--818, 2003.

\bibitem[CLO97]{cox}
David Cox, John Little, and Donal O'Shea.
\newblock {\em Ideals, varieties, and algorithms}.
\newblock Undergraduate Texts in Mathematics. Springer-Verlag, New York, second
  edition, 1997.
\newblock An introduction to computational algebraic geometry and commutative
  algebra.

\bibitem[CLO05]{using}
David~A. Cox, John Little, and Donal O'Shea.
\newblock {\em Using algebraic geometry}, volume 185 of {\em Graduate Texts in
  Mathematics}.
\newblock Springer, New York, second edition, 2005.

\bibitem[Led08]{jpaa}
M.~Lederer.
\newblock The vanishing ideal of a finite set of closed points in affine space.
\newblock {\em J. Pure Appl. Algebra}, 2008.
\newblock To appear.

\bibitem[MB82]{buchmoell}
H.~M. M{\"o}ller and B.~Buchberger.
\newblock The construction of multivariate polynomials with preassigned zeros.
\newblock In {\em Computer algebra (Marseille, 1982)}, volume 144 of {\em
  Lecture Notes in Comput. Sci.}, pages 24--31. Springer, Berlin, 1982.

\bibitem[MMM93]{mmm}
M.~G. Marinari, H.~M. M{\"o}ller, and T.~Mora.
\newblock Gr\"obner bases of ideals defined by functionals with an application
  to ideals of projective points.
\newblock {\em Appl. Algebra Engrg. Comm. Comput.}, 4(2):103--145, 1993.

\bibitem[MR88]{mora}
Teo Mora and Lorenzo Robbiano.
\newblock The {G}r\"obner fan of an ideal.
\newblock {\em J. Symbolic Comput.}, 6(2-3):183--208, 1988.
\newblock Computational aspects of commutative algebra.

\bibitem[Stu96]{sturmfels}
Bernd Sturmfels.
\newblock {\em Gr\"obner bases and convex polytopes}, volume~8 of {\em
  University Lecture Series}.
\newblock American Mathematical Society, Providence, RI, 1996.

\bibitem[Wei89]{weis}
Volker Weispfenning.
\newblock Constructing universal {G}r\"obner bases.
\newblock In {\em Applied algebra, algebraic algorithms and error-correcting
  codes (Menorca, 1987)}, volume 356 of {\em Lecture Notes in Comput. Sci.},
  pages 408--417. Springer, Berlin, 1989.

\end{thebibliography}
\bibliographystyle{alpha}

\end{document}